\documentclass[12pt]{article}

\usepackage{amsfonts,amsmath}
\usepackage{amssymb,latexsym}
\usepackage[dvips]{epsfig}
\usepackage{amscd}  
\usepackage{pstricks,pst-node}

\title{\bf Crossingless matchings and the cohomology of $(n,n)$ 
                 Springer varieties}
\author{Mikhail Khovanov} 
\date{February 12, 2002}

\newtheorem{prop}{Proposition}
\newtheorem{theorem}{Theorem}
\newtheorem{lemma}{Lemma}
\newtheorem{corollary}{Corollary}

\newtheorem{conjecture}{Conjecture} 
\newcommand{\oplusop}[1]{{\mathop{\oplus}\limits_{#1}}}

\begin{document}
\maketitle
\baselineskip 14pt 

\tableofcontents

\def\C{\mathbb C}
\def\R{\mathbb R}
\def\N{\mathbb N}
\def\Z{\mathbb Z}
\def\Q{\mathbb Q}
\def\F{\mathbb F}
\def\H{\mathbb H} 
\def\P{\mathbb P}
\def\l{\lbrace}
\def\r{\rbrace}
\def\o{\otimes}
\def\D{\Delta}
\def\lra{\longrightarrow}
\def\slt{\mathfrak{sl}_2}

\newcommand{\sumop}[1]{{\mathop{\sum}\limits_{#1}}}
\newcommand{\bigoplusop}[1]{{\mathop{\bigoplus}\limits_{#1}}}
\newcommand{\bigcupop}[1]{{\mathop{\bigcup}\limits_{#1}}}
\newcommand{\cupop}[1]{{\mathop{\cup}\limits_{#1}}}

\def\Hom{\mathrm{Hom}}

\def\sbinom#1#2{\left( \hspace{-0.06in}\begin{array}{c} #1 \\ #2 \end{array}
 \hspace{-0.06in} \right)}

\def\mc{\mathcal} 
\def\mf{\mathfrak}

\def\sl{\mathfrak{sl}}

\newcommand{\vsp}{\vspace{0.1in}}

\def\yesnocases#1#2#3#4{\left\{
\begin{array}{ll} #1 & #2 \\ #3 & #4 
\end{array} \right. }

\def\drawing#1{\begin{center} \epsfig{file=#1} \end{center}}
\def\hsm{\hspace{0.05in}}
\newcommand{\cC}{\mathcal{C}}
\newcommand{\Id}{\mathrm{Id}}
\newcommand{\cA}{\mc{A}}

\newcommand{\define}{\stackrel{\mbox{\scriptsize{def}}}{=}}
\newcommand{\spring}{\mc{B}_{n,n}}  
\newcommand{\Stotal}{\widetilde{S}}  
\newcommand{\Ttotal}{\widetilde{T}}  
\newcommand{\cohom}{\mathrm{H}} 
\newcommand{\hsmall}{\hspace{0.05in}} 

\section{Introduction}

We constructed in \cite{me:tangles} a family of rings $H^n, n\ge 0,$ a new 
invariant of tangles, and a conjectural invariant of tangle cobordisms. The 
invariant of a tangle is a complex of $(H^n,H^m)$-bimodules 
(up to chain homotopy equivalences). 
This paper relates $H^n$ to the Springer variety $\spring$ 
of complete flags in $\C^{2n}$ stabilized by a fixed nilpotent operator 
with two Jordan blocks of size $n.$

\begin{theorem} \label{center-spring} The center of $H^n$ is isomorphic to 
the cohomology ring of $\spring$:   
\begin{equation*}
 Z(H^n)\cong \mathrm{H}^{\ast}(\spring,\Z). 
\end{equation*}
\end{theorem} 
Both rings have natural gradings and the isomorphism is grading-preserving. 

Theorem~\ref{center-spring} is proved in a roundabout way, by finding 
generators and defining relations for both rings. Cohomology ring of 
the Springer variety $\mc{B}_{\lambda},$ for a partition $\lambda$ of 
$n,$ is well-understood. In Section~\ref{spring-gen-rel} we use a 
presentation of $\mc{B}_{\lambda}$ via generators and relations obtained 
by de Concini and Procesi \cite{ConciniProcesi} to prove

\begin{theorem}\label{springer-present} 
The cohomology ring of $\mc{B}_{n,n}$ is isomorphic, as a 
graded ring, to the quotient of the 
polynomial ring $R=\Z[X_1, \dots, X_{2n}], \mathrm{deg} X_i =2,$ by the ideal 
$R_1$ with generators 
 \begin{eqnarray}
  \label{sq-zero}           X_i^2,  &  \hspace{0.1in} & i\in [1,2n]; \\
  \label{el-sym} \sum_{|I|=k} X_I,  &  \hspace{0.1in} & k\in [1,2n];
 \end{eqnarray}
where $X_I = X_{i_1}\dots X_{i_k}$ for $I=\{i_1, \dots, i_k\}$ and 
the sum is over all cardinality $k$ subsets of $[1,2n].$  
\end{theorem}

Generators (\ref{el-sym}) are elementary symmetric polynomials in 
$X_1, \dots ,X_{2n}.$ The quotient of $R$ by the ideal generated by 
(\ref{el-sym}) is isomorphic to the cohomology ring of the variety $\mc{B}$ of 
complete flags in $\C^{2n}.$ The inclusion $\spring\subset \mc{B}$ induces 
a surjection of cohomology rings $\cohom^{\ast}(\mc{B},\Z)\longrightarrow 
\cohom^{\ast}(\spring, \Z).$ It turns out that by adding relations 
$X_i^2=0$ we get a presentation for the cohomology ring of $\spring.$ 

\vspace{0.07in}
 
We recall several notations and definitions from \cite{me:tangles}, including 
that of $H^n.$ Let $\mc{A} = \cohom^{\ast}(S^2,\Z)$ be the cohomology ring of 
the 2-sphere, $\mc{A}\cong \Z[X]/(X^2).$ The trace form 
\begin{equation*}
  \mathrm{tr}:\mc{A}\to \Z, \hspace{0.2in} \mathrm{tr}(1)=0, \hspace{0.1in}
  \mathrm{tr}(X)=1,
\end{equation*}
makes $\mc{A}$ into a commutative Frobenius algebra. We assign to $\mc{A}$ a 
2-dimensional topological quantum field theory $\mc{F},$ a functor 
from the category of oriented cobordisms between 1-manifolds to the 
category of abelian groups. $\mc{F}$ associates 
\begin{itemize}
\item $\mc{A}^{\otimes k}$ to the disjoint union of $k$ circles,
\item the multiplication map $\mc{A}^{\otimes 2}\to \mc{A}$ to the "pants" 
   cobordism (three-holed sphere viewed as a cobordism from two circles to 
   one circle),
\item the comultiplication 
 \begin{equation*} 
  \Delta: \mc{A}\to \mc{A}^{\otimes 2}, \hspace{0.2in}
  \Delta(1)= 1\otimes X + X\otimes 1, \hspace{0.1in}
  \Delta(X)= X \otimes X 
 \end{equation*}
to the "inverted pants" cobordism,
\item either the trace or the unit map to the disk (depending on whether 
 we consider the disk as a cobordism from one circle to the empty manifold 
 or vice versa). 
\end{itemize}

\begin{figure} [htb] \drawing{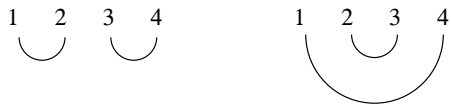}\caption{crossingless matchings 
   $\{(12),(34)\}$   and $\{ (14), (23) \}$} \label{pic-match2} 
 \end{figure}

Let $B^n$ be the set of crossingless matchings of $2n$ points. Equivalently, 
$B^n$ is all pairings of integers from $1$ to $2n$ such that there 
is no quadruple $i<j<k<l$ with $(i,k)$ and $(j,l)$ paired. Most of the 
time $n$ is fixed, and we denote $B^n$ simply by $B.$ Figure~\ref{pic-match2} 
depicts elements of $B^2.$ 

\begin{figure} [htb] \drawing{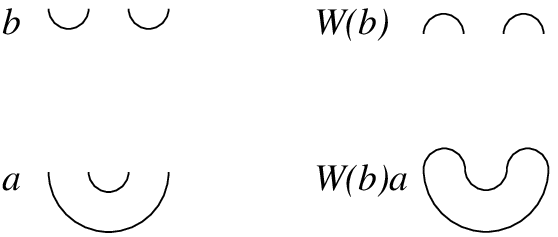}\caption{} \label{pic-glue} 
 \end{figure}

For $a,b\in B$ denote by $W(b)$ the reflection of $b$ about the horizontal 
axis, and by $W(b)a$ the closed 1-manifold obtained by gluing $W(b)$ 
and $a$ along their boundaries, see figure~\ref{pic-glue}. 

\begin{figure} [htb] \drawing{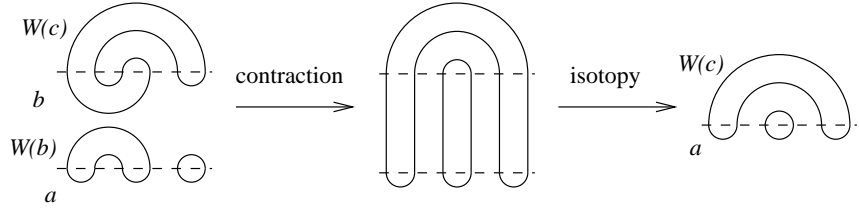}\caption{The contraction cobordism} 
 \label{pic-contr} 
 \end{figure}

$\mc{F}(W(b)a)$ is an abelian group isomorphic to $\mc{A}^{\otimes I},$ 
where $I$ is the set of connected components of $W(b)a.$ 
For $a,b,c\in B$ there is a canonical cobordism from $W(c)bW(b)a$ to 
$W(c)a$ given by "contracting" $b$ with $W(b),$ see figure~\ref{pic-contr}. 

This cobordism induces a homomorphism of abelian groups 
 \begin{equation}\label{ind-hom}
   \mc{F}(W(c)b)\otimes \mc{F}(W(b)a) \lra \mc{F}(W(c)a).
 \end{equation}
Let
 \begin{equation*}
   H^n\define \oplusop{a,b\in B} \hspace{0.05in} {_b(H^n)_a}, 
  \hspace{0.2in} {_b(H^n)_a}\define \mc{F}(W(b)a).  
  \end{equation*} 
Homomorphisms (\ref{ind-hom}), over all $a,b,c,$ define an associative 
multiplication in $H^n$ (we let the products 
 ${_d(H^n)_c}\otimes\hspace{0.05in} {_b(H^n)_a}\to\hspace{0.05in}{_d(H^n)_a}$ 
 be zero if $b\not= c$). 

$_a(H^n)_a$ is a subring of $H^n,$ isomorphic to $\mc{A}^{\otimes n}.$ 
Its element $1_a\define 1^{\otimes n}\in \mc{A}^{\otimes n}$ is an 
idempotent in $H^n.$ The sum $\sum_a 1_a$ is the unit element of $H^n.$ 
Notice that $_b(H^n)_a= 1_b H^n 1_a.$ 

Cohomological grading of $\mc{A}$ ($\mbox{deg}(1)=0, \mbox{deg}(X)=2$) 
  gives rise to a grading of $H^n,$ see \cite{me:tangles} 
for details. 

\vspace{0.1in} 

Denote by $S$ the 2-sphere $\mathbb{S}^2.$ Consider the direct product 
 \begin{equation*}
   S^{\times 2n} \define S\times S \times \dots \times S \hspace{0.2in}
   \mbox{(2n terms)}
 \end{equation*}
and a submanifold $S_a\in S^{\times 2n},$ for $a\in B,$ which consists of 
sequences \newline 
$(x_1, \dots , x_{2n}), x_i\in S$ such that $x_i=x_j$ 
whenever $(i,j)$ is a pair in $a.$ This submanifold is diffeomorphic to 
$S^{\times n}.$ Let $\Stotal = \cup_{a\in B} S_a,$ a subspace 
in $S^{\times 2n}.$ For example, if $n=2$ then $\Stotal$ is homeomorphic to 
two copies of $S\times S$ glued together along their diagonals. 

$H^n$ and $\Stotal$ are constructed along similar lines:  
the cohomology ring of $S_a$ is canonically isomorphic to the ring 
$_a(H^n)_a,$  while the abelian group $_a(H^n)_b$ is canonically isomorphic to 
$\cohom^{\ast}(S_a\cap S_b, \Z).$ These isomorphisms are our starting point 
in the proof (Section~\ref{proof-3}) of  

\begin{theorem} \label{center-config} 
The center of $H^n$ is isomorphic to the cohomology ring of $\Stotal$: 
\begin{equation*}
  Z(H^n) \cong \mathrm{H}^{\ast} (\Stotal, \Z). 
\end{equation*} 
\end{theorem}

In Section~\ref{ZH-g-r} we prove 

\begin{theorem}\label{gen-and-rel} 
Cohomology ring of $\Stotal$ is isomorphic, as a graded ring, to the quotient 
of the polynomial ring $R=\Z[X_1, \dots, X_{2n}], \mathrm{deg} X_i =2,$ 
by relations
 \begin{eqnarray}
  \label{sq-zero2}           X_i^2 & = & 0, \hspace{0.2in} i\in [1,2n]; \\
  \label{el-sym2} \sum_{|I|=k} X_I & = & 0, \hspace{0.2in} k\in [1,2n].
 \end{eqnarray}
\end{theorem}

Theorems~\ref{springer-present}, \ref{center-config} and 
\ref{gen-and-rel} imply Theorem~\ref{center-spring}.  They also show that
spaces $\Stotal$ and $\spring$ have isomorphic cohomology rings. 
These spaces have similar combinatorial structure. Irreducible components 
of $\spring,$ just like those of $\Stotal,$ are enumerated by crossingless 
matchings. Each component $K_a\subset \mc{B}_{n,n}$ is an 
iterated $\mathbb{P}^1$-bundle over $\mathbb{P}^1$ (see \cite{Fung}),  
and homeomorphic to $S_a\subset \Stotal.$ Moreover, $K_a\cap K_b$ and 
  $S_a\cap S_b$ are homeomorphic. We expect that there is a compatible family 
of homeomorphisms and suggest  

\begin{conjecture} $\spring$ and $\Stotal$ are homeomorphic. 
\end{conjecture}

\emph{Warning:}
$\Stotal$ can be naively upgraded to an algebraic variety, by changing the 
2-sphere $S$ to $\mathbb{P}^1$ in the definition of $\Stotal.$ With this 
structure, however, $\Stotal$ 
is not isomorphic to $\spring$ as an algebraic variety, since irreducible 
components of $\spring$ are nontrivial iterated $\mathbb{P}^1$ bundles over 
$\mathbb{P}^1,$ while those of $\Stotal$ are just direct 
products of $\mathbb{P}^1.$ 

\vspace{0.1in}

Let $\mc{B}_{\kappa}$ be the Springer variety of complete flags in $\C^n$ 
stabilized by a fixed nilpotent operator with Jordan decomposition 
$\kappa=(k_1, \dots, k_m).$ 
The cohomology ring of $\mc{B}_{\kappa}$ admits a natural action 
of the symmetric group, see \cite[Section 3.6]{ChrissGinzburg} and 
references therein.
In particular, $S_{2n}$ acts on the cohomology ring of $\spring.$ In view 
of Theorem~\ref{center-spring} it therefore acts on the center 
of $H^n.$ Explicitly, the action is by permutations of $X_i$'s. 
It does not come from any action of $S_{2n}$ on $H^n.$ 

In Section~\ref{symm-group-action} we present an intrinsic construction of 
this action. The $2n$-stranded braid group acts on the category 
$\mc{K}$ of complexes of $H^n$-modules modulo chain homotopies, as 
follows from \cite{me:tangles}. This action 
descends to the braid group action on centers of $\mc{K}$ and 
$H^n$ and factors through to the symmetric group action on the center 
of $H^n.$ 

In Section~\ref{conj-highest-weight} we discuss conjectural 
isomorphisms between centers of 
parabolic blocks of the highest weight category for the Lie algebra
$\mf{sl}_n$ and cohomology algebras of Springer varieties. 

\vspace{0.1in}

\emph{Acknowledgments:} Ideas relating rings $H^n$ and $(n,n)$ 
Springer 
varieties appeared during discussions between Paul Seidel and the author, 
and we are planning a joint paper on various aspects of this correspondence 
\cite{me:Paul2}. The present work can be viewed as a side result 
of \cite{me:Paul2}. I am grateful to Ragnar-Olaf Buchweitz and 
Ivan Mirkovic for useful consultations. This work was partially 
supported by NSF grant DMS-0104139.

%
%

\section{Proof of Theorem~\ref{springer-present}} \label{spring-gen-rel}

De Concini and Procesi \cite{ConciniProcesi} found a presentation for 
the cohomology ring of the Springer variety associated to a partition.
We describe their result specialized to the $(n,n)$ partition. 

Start with the ring $R=\Z[X_1, \dots, X_{2n}].$ 
For $I\subset [1,2n]$ let $X_I=\prod_{i\in I}X_i$ and let
$e_k(I)$ be the elementary symmetric polynomial of order $k$ in variables 
$X_i,$ for $i\in I$: 
 \begin{equation*}
  e_k(I) = \sum_{|J|=k, J\subset I} X_J. 
 \end{equation*}

 \begin{prop} (see \cite{ConciniProcesi}) The cohomology ring of $\spring$ is 
  isomorphic to the quotient ring of $R=\Z[X_1,\dots, X_{2n}]$ by the ideal 
  $R_2$ generated by $e_k(I)$ for all $k+|I|=2n+1,$  $X_I$ for all $|I|=n+1,$ 
  and $X_i^2$ for $i\in [1,2n].$ 
 \end{prop} 

\emph{Remark:} Defining relations are expressed in \cite{ConciniProcesi}
in terms of complete symmetric functions. Complete and elementary symmetric 
functions coincide modulo the ideal generated by $X_i^2, i\in [1,2n].$ 

We want to show the equality of ideals $R_1=R_2,$ where $R_1$ was defined 
in Theorem~\ref{springer-present}. Let $R_3$ be the ideal of $R$ generated by  
$X_i^2, i\in [1,2n].$ Let $e_k=e_k([1,2n])\in R_1.$ 

\emph{Claim:} $R_2\subset R_1.$ Since both $R_1$ and $R_2$ are stable under 
the permutation action of $S_{2n}$ on $R,$ it suffices to show that 
$e_k([1,2n-k+1])$ and $X_{[1,n+1]}$ lie in  $R_1.$ Indeed, 
 \begin{eqnarray}
  \label{equality4} e_k([1,2n-k+1])& \equiv & \sum_{i=0}^{k-1} (-1)^i e_i
  ([2n-k+2, 2n])e_{k-i}\hspace{0.2in}(\mathrm{mod }\hspace{0.1in}R_3), \\
   X_{[1,n+1]} & \equiv & \sum_{i=0}^{n-1} (-1)^i e_i([n+2,2n]) e_{n+1-i} 
  \hspace{0.2in}(\mathrm{mod }\hspace{0.1in}R_3).
  \end{eqnarray} 
The right hand sides lie in $R_1/R_3$ and the claim follows. $\square$ 

\emph{Claim:} $R_1\subset R_2.$ Equalities (\ref{equality4}) and induction 
on $k$ imply that $e_k\in R_2$ for all $k\le n.$ Moreover, $e_k\in R_2$ 
for $k>n$ since $X_I\in R_2$ for any $|I|>n.$ $\square$ 

Therefore, $R_1=R_2$ and the cohomology ring of $\spring$ is isomorphic 
to $R/R_1.$ 

$\square$

%
%

\section{Proof of Theorem~\ref{center-config}}
\label{proof-3}

$\cohom(Y)$ will denote the cohomology ring of a topological space $Y$ 
with integer coefficients. 

There is a canonical isomorphism of rings 
$\cohom(S_a) \cong \cA^{\otimes I}\cong \hspace{0.03in} _aH_a,$ where $I$ is 
the set of arcs of $a.$ Similarly, there are natural abelian group 
isomorphisms $\cohom(S_a\cap S_b) \cong \cA^{\otimes I}\cong 
  \hspace{0.03in} _a H_b,$
where $I$ is the set of connected components of $W(a)b.$   These isomorphisms 
allow us to make $_a H_b$ into a ring. $_a 1_b \define 1^{\otimes k} \in 
\mc{A}^{\otimes k}$ is the unit of this ring.

Inclusions $ S_b\supset (S_a\cap S_b) \subset S_a$ induce ring homomorphisms 
 \begin{equation*}
  \psi_{a;a,b}: \cohom(S_a) \longrightarrow \cohom(S_a\cap S_b),
   \hspace{0.2in} 
  \psi_{b;a,b}: \cohom(S_b) \longrightarrow \cohom(S_a\cap S_b). 
 \end{equation*}
Likewise, maps 
 \begin{equation*}
 \gamma_{a;a,b}:\hspace{0.06in} _a(H^n)_a \longrightarrow \hsmall _a(H^n)_b, 
  \hspace{0.2in}
 \gamma_{b;a,b}:\hspace{0.06in} _b(H^n)_b \longrightarrow  \hsmall _a(H^n)_b, 
 \end{equation*}
  given by $x \longmapsto \hsmall x\hspace{0.03in} _a 1_b $ and 
  $x \longmapsto \hspace{0.03in} _a 1_b x$ are ring homomorphisms.  
The following diagram made out of these homomorphisms and automorphisms 
 commutes: 
 \begin{equation} \label{psi-diag}
    \begin{CD}   
     \cohom(S_a) @>>>  \cohom(S_a\cap S_b) @<<<
     \cohom(S_b)    \\ 
     @VV{\cong}V                @VV{\cong}V   @VV{\cong}V   \\
     _a(H^n)_a   @>>>  _a(H^n)_b           
      @<<< _b(H^n)_b  
    \end{CD} 
 \end{equation}

Suppose given finite sets $I$ and $J,$ rings $A_i,i\in I$ and $B_j, j\in J,$ 
and ring homomorphisms $\beta_{i,j}: A_i\to B_j$ for some pairs $(i,j)$. 
Let 
\begin{equation*}
  \beta= \sum \beta_{i,j}, \hspace{0.2in} \beta: \prod_{i\in I} A_i \lra 
  \prod_{j\in J} B_j.  
\end{equation*}
Define the \emph{equalizer} of $\beta$ (denoted $\mathrm{Eq}(\beta)$) 
as the subring of ${\mathop{\prod}\limits_{i\in I}A_i}$ which consist 
of ${\mathop{\times}\limits_{i\in I}a_i}$ such that 
$\beta_{i,j} a_i = \beta_{k,j}a_k$ whenever 
$\psi_{i,j}$ and $\psi_{k,j}$ are defined. 

\vspace{0.1in}

Diagrams (\ref{psi-diag}) give rise to a commutative diagram 
of ring homomorphisms

 \begin{equation} \label{larger-diag}
    \begin{CD}   
      \mathrm{Eq}(\psi)    @>>>  {\mathop{\prod}\limits_{a}} \hsmall 
     \cohom(S_a) @>{\psi}>>
    {\mathop{\prod}\limits_{a\not= b}} \cohom(S_a\cap S_b) \\ 
      @VV{\cong}V     @VV{\cong}V                @VV{\cong}V     \\
     \mathrm{Eq}(\gamma) @>>> {\mathop{\prod}\limits_{a}} 
     \hsmall _a(H^n)_a @>{\gamma}>> 
    {\mathop{\prod}\limits_{a\not= b}} \hsmall _a(H^n)_b     
    \end{CD} 
   \end{equation} 
where 
\begin{equation*}
\psi = \sum_{a\not= b}(\psi_{a;a,b}+ \psi_{b;a,b})\hspace{0.1in}
 \mathrm{and} \hspace{0.1in}
 \gamma = \sum_{a\not= b}(\gamma_{a;a,b}+ \gamma_{b;a,b}).
\end{equation*}

For an element $z\in H$ write $z= \sum_{a,b}\hsmall _b z_a$ where $_b z_a\in 
\hspace{0.03in} _b H_a.$ If $z$ is central, $_b z_a=0$ if $a\not= b,$ since 
$0 = z 1_b 1_a = 1_b z 1_a =\hsmall _b z_a.$ Thus, 
$z= \sum_a \hsmall _a z_a. $ Denote $_az_a$ by $z_a.$ 
Clearly, $z= \sum_a z_a$ is central iff 
$z_a \hspace{0.02in} _a1_b = 
 \hsmall _a1_b z_b$ for all $a,b$ such that $a\not= b.$  
 Therefore, $Z(H^n)\cong \mathrm{Eq}(\gamma).$ 

\vspace{0.1in}

Inclusions $S_a\subset \Stotal$ induce ring homomorphisms 
$\cohom(\Stotal) \lra {\mathop{\prod}\limits_{a}} S_a$ which factor 
through $\mathrm{Eq}(\psi).$ Putting everything together, we obtain the 
following diagram

 \begin{figure}[htb]
 \begin{pspicture}(0,0)(2.5,3)
 \begin{psmatrix}
   $\cohom(\Stotal)\hsmall$  &  \hsmall $\mathrm{Eq}(\psi)$ \hsmall &  
   \hsmall  ${\mathop{\prod}\limits_{a}} \hsmall \cohom(S_a)\hsmall $ &   
   \hsmall ${\mathop{\prod}\limits_{a\not= b}} \cohom(S_a\cap S_b)$ \hsmall \\
   \hsmall  $Z(H^n)$ \hsmall  & \hsmall $\mathrm{Eq}(\gamma)$ \hsmall &  
   \hsmall  ${\mathop{\prod}\limits_{a}}^{\mathrm{}} 
    \hsmall _a(H^n)_a$ \hsmall  &  
   \hsmall ${\mathop{\prod}\limits_{a\not= b}} \hsmall _a(H^n)_b$ \hsmall  
   \psset{arrows=->}
   \ncline{1,1}{1,2}^{$\tau$} 
   \ncline{1,2}{1,3}
   \ncline{1,3}{1,4}^{$\psi$}
   \ncline{2,1}{2,2}^{$\cong$} 
   \ncline{2,2}{2,3}
   \ncline{2,3}{2,4}^{$\gamma$}
   \ncline{1,2}{2,2}>{$\cong$}
   \ncline{1,3}{2,3}>{$\cong$}
   \ncline{1,4}{2,4}>{$\cong$}
 \end{psmatrix}
 \end{pspicture} 
 \end{figure}

\vspace{0.2in}

Theorem~\ref{center-config} will follow from 

\begin{prop} \label{ta-isom} $\tau$ is an isomorphism. 
\end{prop}

Proof of this proposition occupies the rest of this section.  

\vspace{0.1in}

\begin{figure}[htb]\drawing{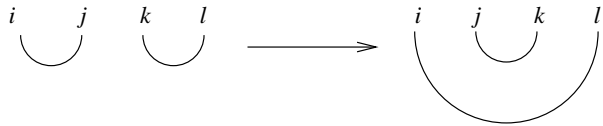} 
  \caption{$a\rightarrow b$} \label{pic-arrow} \end{figure}

\begin{figure}[htb]\drawing{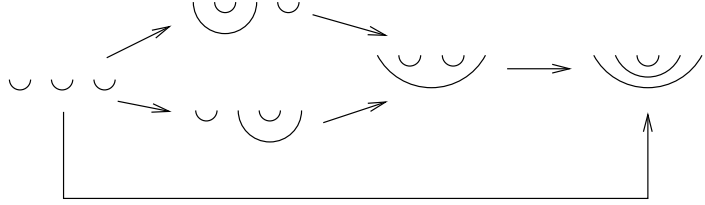}\caption{Arrow relations for $n=3$}
 \label{pic-allarrow} \end{figure}

For $a,b\in B$ we will write $a\rightarrow b$ 
if there is a quadruple $i<j<k<l$ such that $(i,j)$ and $(k,l)$ are 
pairs in $a,$ $(i,l)$ and $(j,k)$ are pairs in $b,$ and otherwise $a$ and 
$b$ are identical (see Figure~\ref{pic-arrow}). 
Figure~\ref{pic-allarrow} depicts all arrow relations for $n=3.$ 

Introduce a partial order on $B$ by $a\prec b$ iff there is 
a chain of arrows $a\rightarrow a_1 \rightarrow \dots \rightarrow a_m 
\rightarrow b.$ 
We extend the partial order $\prec$ to a total order on $B$ in an 
arbitrary way and denote it by $<.$

\vspace{0.2in}

Define the distance $d(a,b)$ between $a$ and $b$ as the minimal length $m$  
of a sequence $(a=a_0,a_1, \dots, a_m=b)$ such that for each $i$ either 
$a_i\rightarrow a_{i+1}$ or $a_{i+1} \rightarrow a_i.$ One geometric 
interpretation of the distance: the diagram $W(b)a$ has 
$n-d(a,b)$ circles. 

\begin{lemma} \label{lemma-sink}
For any $a,b\in B$ there is $c$ such that 
$d(a,b)= d(a,c) + d(c,b)$ and $a \succ  c \prec b.$  
\end{lemma}

Proof is left to the reader. $\square$

\begin{lemma} If $d(a,c) = d(a,b) + d(b,c)$ then 
  \begin{equation*} S_a \cap S_c = S_a \cap S_b \cap S_c. \end{equation*}
\end{lemma}

$\square$ 

\vspace{0.1in}

Let $S_{<a} = \bigcupop{b<a} S_b$ and $S_{\le a} = \bigcupop{b\le a} S_b.$ 
Note that if $c$ is the next element after $a$ in the total order $<$ on $B$
then $S_{<c} = S_{\le a}.$ 

\begin{lemma} 
\begin{equation*}
S_{<a}\cap S_a = \cup_{b\rightarrow a} (S_b\cap S_a).
\end{equation*}
\end{lemma}

Follows from the previous lemma and lemma~\ref{lemma-sink}.  $\square$

\begin{lemma} \label{even-only} $S_{<a}\cap S_a$ has cohomology in even 
degrees only. The inclusion $(S_{<a}\cap S_a) \subset S_a$ induces a 
\emph{surjective} homomorphism of cohomology rings 
$\cohom(S_a)\lra \cohom(S_{<a}\cap S_a).$
\end{lemma} 

\begin{figure} [htb] \drawing{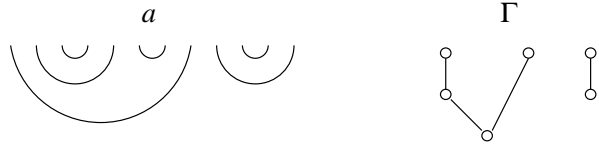} \caption{A crossingless 
 matching and its graph}  
 \label{pic-graph} \end{figure}

\emph{Proof:} We construct a cell decomposition of $S_a.$ Let $I$ be the 
set of arcs of $a.$ There is a canonical homeomorphism 
 $S_a\cong S^{\times I}.$ 
Let $\Gamma$ be the graph with $I$ as the set of vertices and $y,z\in I$ 
are connected by an edge iff there exist $b\rightarrow a$ such that $b$ is 
obtained from $a$ by erazing $y,z$ and reconnecting their endpoints in 
a different way. See Figure~\ref{pic-graph} for an example. 

$\Gamma$ is a forest (a disjoint union of trees). Let $E$ be the set of 
edges of $\Gamma.$ Mark a vertex in each connected component of $\Gamma$ 
and denote by $M$ the set of marked vertices. Note that $|E|+|M|=n.$ 

Fix a point $p\in S.$ For each $J\subset (E\sqcup M)$ let $c(J)$ be the 
subset of $S^{\times I}$ consisting of points $\{ x_i\}_{i\in I}, x_i\in S$ 
such that 
 \[ 
 \begin{array}{ll}
  x_i = x_j       & \mbox{if $(i,j)\in J$}, \\
  x_i \not= x_j   & \mbox{if $(i,j)\notin J$}, \\
  x_i = p         & \mbox{if $i\in M\cap J$}, \\
  x_i \not= p     & \mbox{if $i\in M, i\notin J$}. 
 \end{array}
 \]
Clearly, $S^{\times I}= \sqcup_{J} c(J)$ and $c(J)$ is homeomorphic to 
$\R^{2(n-|J|)}.$ We obtain a decomposition of $S_a\cong S^{\times I}$ 
into even dimensional cells. It restricts to a cell decomposition of 
$S_{<a}\cap S_a,$ the latter a union of cells $c(J)$ such 
that $J\cap E \not= \emptyset.$ The lemma follows, since these 
decompositions give us cochain complexes with zero differentials that 
describe cohomology groups of $S_a$ and $S_{<a}\cap S_a.$ $\square$ 

\begin{lemma} \label{lemma-injective} Homomorphism 
 \begin{equation*}
   \cohom(S_{< a} \cap S_a) \lra \oplusop{b<a} \cohom(S_b\cap S_a) 
 \end{equation*}
induced by inclusions $(S_b\cap S_a)\subset (S_{<a}\cap S_a)$ is injective. 
\end{lemma}

\emph{Proof:} It suffices to check that 
 \begin{equation*}
   \cohom(S_{< a} \cap S_a) \lra \oplusop{b\rightarrow a} \cohom(S_b\cap S_a) 
 \end{equation*}
is injective. The cell decomposition of $S_{<a}\cap S_a$ constructed above 
restricts to a cell decomposition of $S_b\cap S_a,$ for each $b\rightarrow a.$
Since $S_{<a}\cap S_a = \cup_{b\rightarrow a} (S_b\cap S_a),$ the lemma 
follows. $\square$

Note that $S_{\le a} = S_{<a} \cup S_a.$ Consider the Mayer-Vietoris 
sequence for $(S_{<a},S_a):$ 
 \begin{equation*}
  \lra \cohom^m(S_{\le a}) \lra \cohom^m(S_{<a})\oplus \cohom^m(S_a) 
  \lra \cohom^m(S_{<a}\cap S_a) \lra 
 \end{equation*}
 
\vspace{0.1in}

\begin{prop} \label{prop-exact} $S_{\le a}$ has cohomology in even degrees
only. The Mayer-Vietoris sequence for $(S_{<a},S_a)$ breaks down into 
short exact sequences
\begin{equation}\label{seq-less}
0 \lra \cohom^{2m}(S_{\le a}) \lra \cohom^{2m}(S_{<a})\oplus \cohom^{2m}(S_a) 
\lra \cohom^{2m}(S_{<a} \cap S_a)\lra 0
\end{equation}
for $0\le m \le n.$  
\end{prop}

\emph{Proof:} Induction on $a$ with respect to the total order $<.$ 
Induction base is obvious. Induction step: let $e$ be the element before $a$ 
relative to $<$ and assume the proposition holds for $e.$ Then spaces 
$S_{<a} = S_{\le e}, S_a,$ and $S_{<a}\cap S_a$ have cohomology 
in even degrees only (the last one by lemma~\ref{even-only}) and the 
Mayer-Vietoris sequence degenerates into exact sequences 
 \begin{eqnarray*}
  & &  0 \lra \cohom^{2m}(S_{\le a}) \lra \cohom^{2m}(S_{<a})\oplus 
   \cohom^{2m}(S_a)\lra \cohom^{2m}(S_{<a} \cap S_a)\lra \\
  & & \lra   \cohom^{2m+1}(S_{\le a}) \lra 0. 
 \end{eqnarray*}
By lemma~\ref{even-only} the map $\cohom(S_a)\lra \cohom(S_{<a}\cap S_a)$ is 
surjective, so that the last term of the sequence is zero. $\square$

\begin{prop} \label{prop-exact2} The following sequence is exact
\begin{equation}\label{another-exact}
 0 \lra \cohom(S_{\le a}) \stackrel{\phi}{\lra} \oplusop{b\le a} \cohom(S_b) 
  \stackrel{\psi^-}{\lra}  
 \oplusop{b< c\le a} \cohom(S_b\cap S_c), 
\end{equation}
\end{prop} 
where $\phi$ is induced by inclusions $S_b\subset S_{\le a},$ while 
\begin{equation*}
\psi^- \define \sum_{b<c\le a} (\psi_{b,c}-\psi_{c,b}), 
\end{equation*}
where 
\begin{equation*}
\psi_{b,c}: \cohom(S_b) \lra \cohom(S_b\cap S_c)
\end{equation*}
is induced by the inclusion $(S_b\cap S_c )\subset S_b.$ 

\emph{Proof:} Induction on $a.$ The induction base, $a$ is minimal relative 
to $<,$  is obvious. 
Induction step: assume $e$ preceeds $a$ relative to $<$ and the claim is 
true for $e.$ 
Lemma~\ref{lemma-injective} allows us to substitute 
$\oplusop{b<a}\cohom(S_b\cap S_a)$ 
for $\cohom(S_{<a}\cap S_a)$ in the sequence (\ref{seq-less}) while 
maintaining exactness everywhere but in the last term. Thus, 
\begin{equation}\label{seq-plus}
0 \lra \cohom(S_{\le a}) \lra \cohom(S_{<a})\oplus \cohom(S_a) 
\lra \oplusop{b<a} \cohom(S_b \cap S_a)
\end{equation}
is exact. Moreover, $S_{<a}=S_{\le e}.$ By induction hypothesis 
\begin{equation*}
 0 \lra \cohom(S_{\le e}) \lra \oplusop{f\le e} \cohom(S_f) 
  {\lra}  
 \oplusop{f < g\le e} \cohom(S_f\cap S_e) 
\end{equation*}
is exact. Substituting in (\ref{seq-plus}), and using standard properties 
of complexes, we conclude that (\ref{another-exact}) is exact. $\square$ 

When $a$ is the maximal element of $B,$ 
Proposition~\ref{prop-exact2} tells us that the sequence
 \begin{equation}
 0 \lra \cohom(\Stotal) \stackrel{\phi}{\lra} \oplusop{b} \cohom(S_b) 
  \stackrel{\psi^-}{\lra}  
 \oplusop{b< c} \cohom(S_b\cap S_c) 
\end{equation}
is exact. This is equivalent to Proposition~\ref{ta-isom}. $\square$

\vspace{0.1in}

\emph{Remark:} A similar method establishes an isomorphism between the quotient
of $H^n$ by its commutant subspace and homology of $\Stotal$ with 
integer coefficients: 
\begin{equation*}
H^n/[H^n,H^n]\cong \cohom_{\ast}(\Stotal, \Z).
\end{equation*}
If $\Lambda$ is a symmetric ring, the center of $\Lambda$ is 
dual to $\Lambda/[\Lambda,\Lambda].$ Ring $H^n$ is symmetric 
\cite[Section 6.7]{me:tangles}.

%
%

\section{Proof of Theorem~\ref{gen-and-rel}} \label{ZH-g-r}

\begin{lemma} \label{free} $\cohom(\Stotal)$ is a free abelian group of rank 
{\scriptsize $\sbinom{2n}{n}.$}  
\end{lemma} 

\emph{Proof}
The cell decomposition of $S_a$ defined in the proof of 
Lemma~\ref{even-only} restricts to a cell decomposition of $S_a\cap S_{<a}.$ 
Hence, a cell \emph{partition} of $\Stotal$ can be obtained starting 
with the cell decomposition of $S_a,$ for the minimal $a\in B,$ and 
then adding the cells of $S_a\setminus S_{<a},$ over all $a$ in 
$B$ following the total order $<.$ Note that this is a 
cell \emph{partition} of $\Stotal,$ not a cell decomposition, since the 
closure of a cell is not, in general, a union of cells. Nevertheless, since 
all cells are even-dimensional and the boundary of each cell has 
codimension 2 relative to the cell, $\cohom(\Stotal)$ is a free 
abelian group with a basis consisting of delta functions of these cells. 

For $a\in B$ let $t(a)$ be the number of "bottom" arcs of $a,$ that is, arcs  
with no arcs below them. $t(a)$ is also the number of connected components 
of the graph $\Gamma,$ defined in the proof of Lemma~\ref{even-only}. 
For instance, Figure~\ref{pic-match2} diagrams have two and one bottom 
arcs (the left diagram has two). Our decomposition of $S_a\setminus S_{<a}$ 
has $2^{t(a)}$ cells. Therefore, the cell partition of $\Stotal$ has 
$\sum_{a\in B} 2^{t(a)}$ cells. It is easy to see that this sum 
equals {\scriptsize $\sbinom{2n}{n}.$} Lemma follows. $\square$

$\square$

Recall that $X$ denotes a generator of $\cohom^2(S).$ The inclusion 
$\iota:\Stotal \subset S^{\times 2n}$ induces a homomorphism 
$\iota^{\ast}: \cohom(S^{\times 2n}) \to \cohom(\Stotal).$ Let 
$\phi_i: S^{\times 2n}\to S$ be the projection on the $i$-th component. 
Consider the composition $\phi_i\circ \iota $ and let 
\begin{equation*}
X_i \define (-1)^i \iota^{\ast} \circ \phi_i^{\ast} (X), \hspace{0.2in}
  X_i \in \cohom(\Stotal). 
\end{equation*}

We denote by $[1,2n]$ the set of integers from $1$ to $2n.$ For 
$I\subset [1,2n]$ let $X_I = \prod_{i\in I} X_i.$

\begin{prop} \label{ring-gen-def} The cohomology ring of $\Stotal$ is 
generated by $X_i, i\in [1,2n]$ and has defining relations 
 \begin{eqnarray}
  \label{sq-zero-s}           X_i^2 & = & 0, \hspace{0.2in} i\in [1,2n]; \\
  \label{el-sym-s} \sum_{|I|=k} X_I & = & 0, \hspace{0.2in} k\in [1,2n].
 \end{eqnarray}
\end{prop}

\emph{Proof:} First we show that these relations hold. (\ref{sq-zero-s}) 
is obvious. Let $j_a: S_a\subset \Stotal$ and $j_a^{\ast}$ be the induced 
map on cohomology. (\ref{el-sym-s}) will follow if we check that 
\begin{equation}\label{eq-for-a} 
 \sum_{|I|=k} j_a^{\ast} (X_I) =0
\end{equation}
for all $a\in B,$ since 
 \begin{equation*}
  \sum_{a\in B} j_a^{\ast}: \cohom(\Stotal) \lra \oplusop{a\in B} 
  \cohom(S_a)
 \end{equation*}
is an inclusion. If $(i,i')$ is a pair in $a$ then $j_a^{\ast}(X_i X_{i'})=0$ 
and $j_a^{\ast}(X_i + X_{i'})=0$ (because of the term $(-1)^i$ in the 
definition of $X_i,$ and since $i+i'\equiv 1 (\mbox{mod }2)$). Therefore, 
\begin{equation*}
  \sum_{|I|=k, \{i,i'\}\cap I \not= \emptyset} j_a^{\ast} (X_I) =0
\end{equation*}
where the sum is over all subsets of cardinality $k$ that intersect 
$\{i,i'\}$ nontrivially. To take care of the remaining terms in the 
L.H.S. of (\ref{eq-for-a}),  
 \begin{equation*}
  \sum_{|I|=k, \{i,i'\}\cap I = \emptyset} j_a^{\ast} (X_I),
 \end{equation*}
pick another pair $(r, r')$ in $a$ and apply the same reduction to 
it. After $\frac{n-k}{2}+1$ iterations all cardinality $k$ subsets will be 
accounted for. (\ref{eq-for-a}) follows. 

\vspace{0.07in}

We say that a subset $I$ of $[1,2n]$ is \emph{admissible} if $I\cap [1,m]$ 
has at most $\frac{m}{2}$ elements for each $m\in [1,2n].$

 \begin{lemma} \label{number-adm}
   There are {\scriptsize $\sbinom{2n}{n}$} admissible subsets. 
 \end{lemma}

\emph{Proof} is left to the reader. $\square$ 

 \begin{lemma} \label{lin-combin} $X_J,$ for any $J\subset [1,2n],$ is a 
 linear combination of $X_I,$ over admissible $I.$  
 \end{lemma}

 \emph{Proof:} let $y(J)= \sum_{j\in J} j.$ Assume the lemma is false, 
 and find such a non-admissible $J$ with the minimal possible $y(J).$ 
 Take the smallest possible $m$ such that $|J\cap [1,m]| > \frac{m}{2}.$
 Then $m$ is odd, $m=2r+1$ and $|J\cap [1,m]|=r+1.$ 
 Arguments in the proof of Theorem~\ref{springer-present} in 
 Section~\ref{spring-gen-rel} imply that 
 $e_{r+1}((J\cap [1,m])\cup [m+1,2n]) =0.$ Therefore, $X_{J\cap [1,m]}$ 
 is a linear combination of $X_K$ with $K\subset (J\cap [1,m])\cup [m+1,2n]$ 
 and $y(K)> y(J).$ Since $X_J= X_{J\cap [1,m]} X_{J\cap [m+1,2n]},$ 
  this contradicts minimality of $y(J).$   $\square$

 \begin{lemma} \label{linear-ind} $\{ X_I\},$ over all admissible $I,$ are 
   linearly independent in $\cohom(\Stotal).$
 \end{lemma}

\emph{Sketch of proof:} Induction on $n,$ use homomorphism  
 $\cohom(\Stotal_n) \lra \cohom(\Stotal_{n-1})$ induced by the inclusion 
 $\Stotal_{n-1}\subset \Stotal_n$ (where $\Stotal_n$ is what we usually 
 call $\Stotal$). Details are left to the reader. $\square$ 

Lemmas~\ref{free}, \ref{number-adm}, \ref{lin-combin} and 
\ref{linear-ind} imply Proposition~\ref{ring-gen-def} and the following 
results.  

 \begin{corollary} $\cohom(\Stotal)$ has a basis $\{ X_I\},$ over all 
  admissible $I.$ 
 \end{corollary}

 \begin{corollary} The inclusion $\Stotal \subset S^{\times 2n}$ induces a 
   \emph{surjective} ring homomorphism 
   $\cohom(S^{\times 2n}) \lra \cohom(\Stotal).$ 
 \end{corollary}

%
%

\section{Symmetric group action on the center of $H^n$} 
\label{symm-group-action}

\vspace{0.1in}

{\bf The center of a category}

\vspace{0.1in}

The center of a category is defined as the commutative monoid of natural 
transformations of the identity functor. If the category is pre-additive 
($\mbox{Hom}(X,Y)$ is an abelian group for any objects $X,Y,$ and the 
composition of morphisms is bilinear) then the center is a commutative ring. A 
down-to-earth example: the center of the category of modules 
over a ring $A$ is isomorphic to the center of $A.$ 

Let $F$ be a functor in a category $\mc{C}.$ The center of $\mc{C}$ 
acts in two ways on the set $\mbox{End}(F)$ of endomorphisms of $F,$ since
we can compose $\alpha\in Z(\mc{C})$ and $\beta\in \mbox{End}(F)$
on the left or on the right: 
\begin{eqnarray*}
 \alpha\circ\beta & : & F \cong \mathrm{Id}\circ F 
  \stackrel{\alpha\circ\beta}{\lra} \mathrm{Id}\circ F \cong F, \\
 \beta\circ\alpha & : & F \cong F\circ\mathrm{Id}  
  \stackrel{\beta\circ\alpha}{\lra} F\circ\mathrm{Id} \cong F. 
\end{eqnarray*}
Assume that $F$ is invertible. Then any endomorphism of $F$ has the form 
$\mbox{Id}_F\circ \alpha,$ for a unique $\alpha\in Z(\mc{C}),$ as well as
$\alpha'\circ \mbox{Id}_F,$ for a unique $\alpha'\in Z(\mc{C}).$   
Thus, $F$ defines an automorphism of the monoid $Z(\mc{C})$ which takes 
$\alpha$ to $\alpha'.$  

Suppose group $G$ acts weakly on $\mc{C},$ meaning that there are functors 
$F_g: \mc{C}\to \mc{C}$ for each $g\in G,$ such that $F_1\cong \mbox{Id}$ and 
$F_gF_h\cong F_{gh}$ (we do not impose compatibility conditions on these 
isomorphisms). Each $F_g$ is invertible and gives rise to an automorphism 
of the 
center of $\mc{C}.$ We get an action of $G$ on $Z(\mc{C}).$ Therefore, a weak 
group action on a category descends to an action on the center of the 
category. 

\vspace{0.1in}

{\bf Centers of triangulated and derived categories} 

\vspace{0.1in}

Define the center of a triangulated category $\mc{D}$
as the set of natural transformations of the identity functor that 
commute with the shift functor $[1].$ 
Let $\mc{C}$ be an abelian category and $\widehat{\mc{C}}$ one of the 
triangulated categories associated to $\mc{C}$ (for instance, the 
bounded derived category of $\mc{C},$ or the category of bounded 
complexes of objects of $\mc{C}$ modulo chain homotopies). 
There are ring homomorphisms 
\begin{equation}\label{center-hom} 
   Z(\mc{C})\stackrel{f}{\lra}Z(\widehat{\mc{C}})\stackrel{g}{\lra}Z(\mc{C})  
\end{equation}
whose composition is the identity. $f$ extends $\alpha\in Z(\mc{C})$ 
termwise to complexes of objects of $\mc{C}.$ Homomorphism $g$ is induced 
by the inclusion of categories $\mc{C} \subset \widehat{\mc{C}}.$ 

$f$ and $g$ are not always isomorphisms, as observed by Jeremy Rickard. 
If $\mc{C}$ is the category of modules over the exterior 
algebra in one generator, then $f,$ respectively $g,$ has a nontrivial 
cokernel, respectively kernel. 

\emph{Remark:} Ragnar-Olaf Buchweitz pointed out to me that a triangulated 
category $\mc{D}$ also has \emph{extended center} (or "Hochschild cohomology"),
\begin{equation*}
 \oplusop{m\in \Z} \mbox{Hom}_{\mathrm{Fun}(\mc{D})}(\mbox{Id},\mbox{Id}[m]), 
\end{equation*}
(only natural transformations that supercommute with the shift functor 
are included). The extended center of the derived category of 
$\Lambda\mbox{-mod},$ for an algebra $\Lambda,$  
contains the Hochschild cohomology algebra $\mbox{Ext}_{\Lambda^{e}}(\Lambda,
 \Lambda).$  

An action of $G$ on $\widehat{\mc{C}}$ induces an action on
$Z(\widehat{\mc{C}}).$ If $\mathrm{ker}(g)$ is $G$-stable, the 
action descends to $Z(\mc{C}).$ 
Let $\Lambda$ be a ring and $D(\Lambda)$ the bounded derived category 
of $\Lambda\mbox{-mod}.$ If a  self-equivalence $F$ of $D(\Lambda)$ 
is given by tensoring with a bounded complex of left and 
right projective $\Lambda$-bimodules, then it descends to an automorphism of 
$Z(\Lambda)$ (compare with \cite[Proposition 9.2]{Rickard3}).

\vspace{0.1in}

{\bf Symmetric group action}

\vspace{0.1in}

\begin{figure} [htb] \drawing{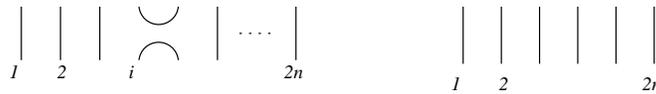}\caption{Diagrams $U_i$ and 
 $Vert_{2n}$} \label{Uid-diag} 
 \end{figure}

The symmetric group action on $Z(H^n)$ can be described intrinsically as 
follows. Let $\mc{K}$ be one of triangulated categories associated to 
 $H^n$ (say, the category of bounded complexes of left 
$H^n$-modules up to chain homotopies).  
The structures described in \cite{me:tangles} lead to a weak action 
of the braid group with $2n$-strands
 on $\mc{K}.$ Diagram $U_i$ (see Figure~\ref{Uid-diag}) defines an 
$H^n$-bimodule $\mc{F}(U_i)$ together with bimodule homomorphisms 
\begin{equation*}
  \mc{F}(U_i) \stackrel{\alpha}{\lra} H^n, \hspace{0.2in}
   H^n \stackrel{\beta}{\lra} \mc{F}(U_i),
\end{equation*} 
induced by elementary cobordisms between $U_i$ and $Vert_{2n}$ 
(we use notations from \cite{me:tangles}, note 
that $H^n\cong \mc{F}(Vert_{2n})$).  

Let $\mc{R}_i: \mc{K}\to \mc{K}$ be the functor of tensoring with 
the complex of bimodules 
\begin{equation} \label{inv-complex}
 0 \lra \mc{F}(U_i) \stackrel{\alpha}{\lra} H^n \lra 0. 
\end{equation} 
$\mc{R}_i$'s are invertible and satisfy the braid group relations. 

\vspace{0.05in}

In Section~\ref{ZH-g-r} we defined generators $X_i$ of $\cohom(\Stotal).$ 
We now describe the image of $X_i$ (also denoted $X_i$) in $Z(H^n)$ under 
the isomorphism $\cohom(\Stotal) \cong Z(H^n)$ established in 
Section~\ref{proof-3}: 
\begin{eqnarray*} 
X_i & = & {\mathop{\sum}\limits_{a\in B}}\hsmall _a(X_i)_a, \hspace{0.2in} 
  _a(X_i)_a \in \hsmall _a(H^n)_a \cong \mc{A}^{\o n}, \\
  _a(X^i)_a & = & (-1)^i 1^{\otimes (n-1)} \otimes X \in \mc{A}^{\otimes (n-1)}
 \otimes \mc{A}, 
\end{eqnarray*}
where the separated $\mc{A}$ corresponds to the circle in $\mc{F}(W(a)a)$ 
that contains the $i$-th endpoint of $a,$ counting from the left. 

For a complex $V$ of $H^n$-bimodules and $z\in Z(H^n)$ let $l_z,$ respectively 
$r_z,$ be the endomorphism of $V$ given by left, respectively right, 
multiplication by $z.$   
Endomorphisms $l_{X_i} - r_{X_{i+1}}$ and $l_{X_{i+1}} - r_{X_i}$ of 
the complex (\ref{inv-complex}) are homotopic to $0,$ via the homotopy 
$\pm \beta: H^n \lra \mc{F}(U_i).$ Therefore, the braid group 
action on $\mc{K}$ descends to the action of the symmetric group 
$S_{2n}$ on $Z(H^n)$ by permutations of $X_i$'s.

%
%

\section{Conjectures on centers of highest weight categories} 
\label{conj-highest-weight}

If $e$ is an idempotent in a ring $\Lambda,$ there is a homomorphism 
$Z(\Lambda)\lra Z(e\Lambda e)$ which takes $z\in Z(\Lambda)$ to $ze.$ 

Let $\mc{O}^{n,n}$ be the full subcategory of a regular block of 
the highest weight category fof $\mf{sl}_{2n}$ which consists of 
locally $U\mf{p}(n,n)$-finite modules, where $\mf{p}(n,n)$ is the 
parabolic subalgebra in $\mf{sl}_{2n}$ of $(n,n)$ block-uppertriangular 
matrices. $\mc{O}^{n,n}$ is equivalent to the category of perverse sheaves 
on the Grassmannian of $n$-planes in $\C^{2n},$ constructible relative 
to the Schubert stratification. There is a unique finite-dimensional 
$\C$-algebra $A_{n,n}$ such that 

 (i) $\mc{O}^{n,n}$ is equivalent to the category of finite-dimensional 
  $A_{n,n}$-modules, 

 (ii) every irreducible $A_{n,n}$-module is one-dimensional. 

$A_{n,n}$ was explicitly described by Tom Braden \cite{Braden}. 
In \cite{me:withTom} we'll construct an idempotent $e$ in $A_{n,n}$ and 
an isomorphism $h: H^n\otimes_{\Z}\C\cong eA_{n,n}e.$ 

 \begin{conjecture} $Z(A_{n,n})\cong Z(eA_{n,n}e),$ and $h$ induces an 
  isomorphism of centers of $H^n\otimes_{\Z}\C$ and $\mc{O}^{n,n}.$ 
 \end{conjecture}

We would like to suggest a more general conjecture relating parabolic highest 
weight categories and Springer varieties. 
Let $\kappa=(k_1, \dots , k_m)$ be a decomposition of $n,$ 
$k_1+ \dots +k_m=n$ and denote by $\mf{p}(\kappa)$ the corresponding 
Lie algebra of block upper-triangular $n$-by-$n$ matrices. Let 
$\mc{O}^{\kappa}$ be the full subcategory of locally 
$U\mf{p}(\kappa)$-finite modules in a regular block $\mc{O}_{reg}$ 
of the highest weight category  $\mc{O}$ for $\mf{sl}_n.$ 
Let $Y^{\kappa}$ be the partial flag variety associated to $\kappa.$ 
It is known that $\mc{O}^{\kappa}$ is equivalent to the category of 
perverse sheaves on $Y^{\kappa},$ smooth along Schubert cells. 

Let $\mc{B}_{\kappa}$ be the Springer variety of complete flags in $\C^n$ 
stabilized by a fixed nilpotent operator with Jordan decomposition 
$(k_1, \dots, k_m).$ 

\begin{conjecture} The center of $\mc{O}^{\kappa}$ is isomorphic 
to the cohomology algebra of $\mc{B}_{\kappa}$: 
\begin{equation}\label{conj-iso1}
Z(\mc{O}^{\kappa}) \cong \mathrm{H}^{\ast}(\mc{B}_{\kappa}, \C).
\end{equation}
\end{conjecture}

Note that the right hand side of the isomorphism (\ref{conj-iso1}) 
depends only on the partition type of $\kappa,$ i.e. preserved by 
permutations of terms 
$k_1, \dots, k_m$ of $\kappa.$ The category $\mc{O}^{\kappa},$ 
featured in the left hand side, generally does not possess the same kind of 
invariance: 

 \begin{prop} \label{inequivalent}
   Categories $\mc{O}^{2,1,1}$ and $\mc{O}^{1,2,1}$ are inequivalent. 
 \end{prop}
 
However, we have

 \begin{prop} \label{der-equiv} If decompositions $\kappa$ and $\kappa'$   
  differ only by a permutation of terms, the categories 
  $\mc{O}^{\kappa}$ and $\mc{O}^{\kappa'}$ are derived equivalent. 
 \end{prop}

Propositions~\ref{inequivalent} and \ref{der-equiv} are proved at the end 
of this section. 
  
If two rings are derived equivalent, their centers
are isomorphic (\cite[Proposition 9.2]{Rickard3}). Since 
$\mc{O}^{\kappa}$ and $\mc{O}^{\kappa'}$ are equivalent to categories 
of modules over finite-dimensional algebras, and these algebras are derived 
equivalent, centers of  $\mc{O}^{\kappa}$ and $\mc{O}^{\kappa'}$ are
isomorphic. Thus, the left hand side of $(\ref{conj-iso1})$ also depends only 
on the partition type of $\kappa.$

We would like conjectural isomorphisms (\ref{conj-iso1}) to be compatible 
with the inclusions of categories $\mc{O}^{\kappa}\subset \mc{O}_{reg}$ and 
topological spaces $\mc{B}_{\kappa}\subset \mc{B},$ where $\mc{B}$ is the 
variety 
of complete flags in $\C^n.$ These inclusion induce ring homomorphisms 
$Z(\mc{O}_{reg})\lra Z(\mc{O}^{\kappa})$ and 
$\cohom^{\ast}(\mc{B},\C)\lra \cohom^{\ast}(\mc{B}_{\kappa},\C)$ which should 
be a part of the following commutative diagram

  \begin{equation} 
    \begin{CD}   
     \cohom^{\ast}(\mc{B},\C) @>>>  \cohom^{\ast}(\mc{B}_{\kappa},\C) \\
     @VV{\cong}V                @VV{\cong}V     \\
     Z(\mc{O}_{reg}) @>>>    Z(\mc{O}^{\kappa})    
    \end{CD} 
   \end{equation} 

\vspace{0.3in}

Let $\Theta_i: \mc{O}^{\kappa}\lra \mc{O}^{\kappa}$ be the functor of 
translation across the $i$-th wall. $\Theta_i$ is the product of 
two biadjoint functors (translations on and off the $i$-th wall). 
Let $\alpha_i : \Theta_i \lra \mathrm{Id}$ be one of the natural 
transformations 
coming from biadjointness. $\mc{R}_i,$ the cone of $\alpha_i,$ is 
a functor in the derived category of $\mc{O}^{\kappa}.$ Functors 
$\mc{R}_i$ are invertible and generate a braid group action in 
$D^b(\mc{O}^{\kappa}).$ We conjecture that this action descends to 
a symmetric group action on $Z(\mc{O}^{\kappa}),$ and the ring isomorphism 
(\ref{conj-iso1}) can be made $S_n$-equivariant.

\vspace{0.1in}

Let $A_{\kappa}$ be a finite-dimensional algebra such that 
$A_{\kappa}\mbox{-mod}\cong \mc{O}^{\kappa}.$ Let $e\in A_{\kappa}$ be 
the maximal idempotent such that the left $A_{\kappa}$-module 
$A_{\kappa}e$ is injective. 

\begin{conjecture} Inclusion $eA_{\kappa}e\subset A_{\kappa}$ induces an 
isomorphism of centers $Z(A_{\kappa}) \cong  Z(e A_{\kappa}e).$ 
\end{conjecture} 

\vspace{0.1in}

\emph{Proof of Proposition~\ref{inequivalent}:} 
Assume the two categories are equivalent. Any intrinsic homological 
information about them is identical. The equivalence restricts to a bijection 
between isomorphism classes of simple objects. The bijection induces 
isomorphisms between Ext rings of simple objects of these categories. 
                        
Simple objects are in a one-to-one correspondence with Schubert cells in 
partial flag varieties $X_{2,1,1}$ and $X_{1,2,1}.$ For a simple object $L$ 
let $IC(L)$ be the intersection cohomology sheaf on the closure of the 
Schubert cell associated to $L.$ Then $\mbox{Ext}(L,L)\cong 
\mbox{Ext}(IC(L),IC(L)).$ 
                        
Let us count the number of simple objects $L$ in each category with 
$\mbox{dim}(\mbox{Ext}(L,L))=3.$ The Schubert cell of such an object is 
necessarily 2-dimensional and its closure is diffeomorphic to 
$\C\mathbb{P}^2$ (use that the cohomology of the closure of the cell 
is a direct summand of $\mbox{Ext}(IC(L),IC(L)).$
$X_{2,1,1}$ has only one such cell, while $X_{1,2,1}$ has two. 
Contradiction. $\square$           
                 
\vspace{0.07in}
                 
\emph{Proof of Proposition~\ref{der-equiv}:}
To construct an equivalence 
between $D^b(\mc{O}^{\kappa})$ and $D^b(\mc{O}^{\kappa'})$
note that these categories are isomorphic to the derived categories 
of sheaves on partial flag varieties $Y_{\kappa}$ and $Y_{\kappa'},$
smooth along Schubert stratifications. It suffices to treat the case when 
$\kappa$ and $\kappa'$ differ by a transposition of adjacent terms, 
$\kappa=(k_1, \dots, k_m), \kappa'=(k_1, \dots, k_{i-1},k_i, \dots, k_m).$  
Let $U\subset Y_{\kappa}\times Y_{\kappa'}$ be the set 
\begin{equation*}
 ((F_1, \dots F_m)\in Y_{\kappa}, (F_1', \dots, F_m')\in Y_{\kappa'})| 
 F_j=F_j' \hspace{0.1in}\mathrm{ for }\hspace{0.1in}
 j\not=i, \hsmall F_i\cap F_i' = F_{i-1},
\end{equation*}
of pairs of partial flags.
Let $\mc{G}$ be sheaf on $Y_{\kappa}\times Y_{\kappa'}$ which is the 
continuation by $0$ of the constant sheaf on $U.$ 
Convolution with $\mc{G}$ is an equivalence of derived categories of sheaves 
on $Y_{\kappa}$ and $Y_{\kappa'},$ and restricts to an equivalence of 
subcategories of cohomologically constructible (relative to the 
Schubert startification) complexes of sheaves. The latter categories are 
equivalent to the derived categories of $\mc{O}^{\kappa}$ and 
$\mc{O}^{\kappa'}.$     $\square$

\bibliographystyle{abbrv} 
\bibliography{$HOME/bibl/me,$HOME/bibl/all}

\end{document}